\documentclass[12pt]{article}
 \usepackage{amsmath}
 \usepackage{amssymb}
\usepackage{amsthm}
\usepackage{graphicx}
\newtheorem{thm}{Theorem}[section]
\newtheorem{cor}[thm]{Corollary}
\newtheorem{lem}[thm]{Lemma}
\newtheorem{prop}[thm]{Proposition}

\newtheorem{rem}[thm]{Remark}

\renewcommand{\a}{\alpha}
\renewcommand{\b}{\beta}

\title{Bounds on Tur{\'a}n determinants  }

\author{Christian Berg\thanks{The present work was done while the first
     author was visiting University of Wroc{\l}aw granted by the HANAP
    project mentioned under the second author.}
\and Ryszard Szwarc
 \thanks{The second author was supported by European Commission Marie
Curie Host
Fellowship for the Transfer of Knowledge ``Harmonic Analysis, Nonlinear
Analysis and
Probability''  MTKD-CT-2004-013389 and by MNiSW Grant N201 054 32/4285.}
}
\begin{document}
\maketitle

\begin{abstract} Let $\mu$ denote a symmetric probability measure on
  $[-1,1]$ and let $(p_n)$ be the corresponding orthogonal polynomials
  normalized such that $p_n(1)=1$. We prove that the
  normalized Tur{\'a}n determinant $\Delta_n(x)/(1-x^2)$, where
  $\Delta_n=p_n^2-p_{n-1}p_{n+1}$, is a Tur{\'a}n determinant of order
  $n-1$ for orthogonal polynomials with respect to $(1-x^2)d\mu(x)$. We
  use this to prove lower and upper  bounds for the normalized
  Tur{\'a}n determinant in the interval $-1<x<1$. 
  \end{abstract}


\section{Introduction}
 In the following we will deal with polynomial sequences $(p_n)$ satisfying

\begin{eqnarray}
&&xp_n(x)=\gamma_np_{n+1}(x)+\alpha_np_{n-1}(x) ,\ n\ge 0,\nonumber\\
&&\alpha_n+\gamma_n=1,\;\alpha_n>0,\;\gamma_n>0,\ n\ge 1,\label{eq:rec} \\
&&p_0(x)=1,\;\alpha_0=0,\ 0<\gamma_0\le 1.\nonumber\end{eqnarray}
Note that $(p_n)$ is uniquely determined by \eqref{eq:rec}
from the recurrence coefficients $\alpha_n,\gamma_n$. It is well-known
that the polynomials $p_n$ are orthogonal with respect to a symmetric
probability measure $\mu$ with compact support.

Define the Tur\'an
determinant by
\begin{equation}\label{eq:Turan}
\Delta_n(x)=p_n^2(x)-p_{n-1}(x)p_{n+1}(x),\;n\ge 1.
\end{equation}

In \cite{Szw} the second author proved non-negativity of the Tur{\'a}n
determinant \eqref{eq:Turan} under certain monotonicity conditions on
the recurrence coefficients, thereby obtaining results for new classes
of polynomials and unifying old results.

If $\gamma_0=1$ the polynomials satisfy $p_n(1)=1$ and therefore the
{\it normalized} Tur{\'a}n determinant $\Delta_n(x)/(1-x^2)$ is a
polynomial in $x$.

We shall prove estimates of the form
\begin{equation}\label{eq:purpose}
c\Delta_n(0)\le \frac{\Delta_n(x)}{1-x^2}\le C\Delta_n(0),\quad -1<x<1,
\end{equation}
under certain regularity conditions on the recurrence coefficients. We
prove e.g. an inequality of the left-hand type if $(\a_n)$ is
increasing and concave, see Theorem \ref{thm:thm2}. In Theorem
\ref{thm:thm2a} we give
an inequality of the right-hand type.

Our results depend on a simple relation between the Tur{\'a}n determinants
of order $n$ and $n-1$ (Proposition \ref{thm:pro1})
and the following observation: The normalized Tur{\'a}n determinant
is essentially a Tur{\'a}n
determinant of order $n-1$ for the polynomials $(q_n)$ defined by
\eqref{eq:q}  below, and if $\mu$ is the orthogonality measure of $(p_n)$,
then $(q_n)$ are orthogonal with respect to
the measure $(1-x^2)d\mu(x)$. See Theorem \ref{thm:thm1} and Remark
\ref{thm:TuranTuran} for a precise statement.

In Proposition  \ref{thm:Turanq} we prove non-negativity of the
Tur{\'a}n determinant for the normalized polynomials $q_n(x)/q_n(1)$
provided the sequence $(\a_n)$ is increasing and concave (or under the
weaker condition \eqref{eq:qconcave}).

Our work is motivated by results about ultraspherical polynomials,
which we describe next.

For $\a>-1$ let $R_n^{(\a,\a)}(x)=P_n^{(\a,\a)}(x)/P_n^{(\a,\a)}(1)$
denote the symmetric Jacobi polynomials normalized to be 1 for $x=1$, i.e.
\begin{equation}\label{eq:Jacobinorm}
R_n^{(\a,\a)}(x)=\frac{(-1)^n}{2^n(\a+1)_n}(1-x^2)^{-\a}\frac{d^n}{dx^n}(1-x^2)^{n+\a},
\end{equation}  
cf. \cite{Sz}. We use the Pochhammer symbol
$(a)_n=a(a+1)\cdot\ldots\cdot(a+n-1)$. The polynomials are orthogonal
with respect to the symmetric weight function
$c_\a(1-x^2)^\a$ on $]-1,1[$. Here
 $1/c_\a=B(\a+1,1/2)$, so the weight is a probability density.
 We have $R_n^{(\a,\a)}(x)=P_n^{(\lambda)}(x)/P_n^{(\lambda)}(1)$ with
 $\a=\lambda-\tfrac12$, where $(P_n^{(\lambda)})$ are the ultraspherical
 polynomials in the notation of \cite{Sz}.

The corresponding Tur{\'a}n determinant of order $n$
\begin{equation}\label{eq:Turan1}
\Delta^{(\a)}_n(x)=R_n^{(\a,\a)}(x)^2-R_{n-1}^{(\a,\a)}(x)R_{n+1}^{(\a,\a)}(x),
\end{equation}
is clearly a polynomial of degree $n$ i $x^2$ and divisible by
$1-x^2$ since it vanishes for $x=\pm 1$. The  following Theorem was
proved in \cite[pp. 381-382]{T:N} and in \cite[sect. 6]{V:L}:

\begin{thm}\label{thm:ultraspherical} The normalized Tur{\'a}n determinant
\begin{equation}\label{eq:Turannorm}
f_n^{(\a)}(x):=\Delta^{(\a)}_n(x)/(1-x^2)
\end{equation}
 is
\begin{enumerate}
\item[(i)] strictly increasing for $0\le x<\infty$ when $\a>-1/2$.
\item[(ii)] constant equal to $1$ for $x\in\mathbb R$ when $\a=-1/2$.
\item[(iii)] strictly decreasing for $0\le x<\infty$ when $-1<\a<-1/2$.
\end{enumerate}
\end{thm}
It is easy to evaluate $f_n^{(\a)}$ at $x=0,1$ giving
\begin{equation}\label{eq:eval01}
f_n^{(\a)}(0)=\mu^{(\a)}_{[n/2]}\mu^{(\a)}_{[(n+1)/2]},\quad 
f_n^{(\a)}(1)=\frac{1}{2\a+2},
\end{equation}
where we have used the notation from \cite{A:G:K:L}
\begin{equation}\label{eq:genmidbinom}
\mu_n^{(\a)}=\frac{\mu_n}{\binom{n+\a}{n}},
\end{equation}
and $\mu_n$ is the normalized binomial mid-coefficient
\begin{equation}\label{eq:midbinom}
\mu_n=2^{-2n}\binom{2n}{n}=\frac{1\cdot 3\cdot
  5\cdot\ldots\cdot(2n-1)}
{2\cdot 4\cdot\ldots\cdot(2n)}.
\end{equation}

\begin{cor}\label{thm:bounds} For $-1<x<1$ we have
\begin{equation}\label{eq:bounds}
f_n^{(\a)}(0)(1-x^2) < \Delta^{(\a)}_n(x) < f_n^{(\a)}(1)(1-x^2) \mbox{
  for } \a>-1/2
\end{equation}
while the inequalities are reversed when $-1<\a<-1/2$. (For $\a=-1/2$
all three terms are equal to $1-x^2$.)
\end{cor}

For $\a=0$ the inequalities \eqref{eq:bounds} reduce to ($-1<x<1$)
\begin{equation}\label{eq:Turan/Le}
\mu_{[n/2]}\mu_{[(n+1)/2]}(1-x^2)< P_n(x)^2-P_{n-1}(x)P_{n+1}(x)< \frac12(1-x^2) 
\end{equation}
for Legendre polynomials $(P_n)$. This result was recently published
in  \cite{A:G:K:L} using a SumCracker Package by Manuel Kauers,
and it was conjectured  that  the monotonicity result remains true for
ultraspherical polynomials when $\a\ge -1/2$. Clearly the authors 
have not been aware of the early results above.\footnote{Motivated by
  this conjecture the present authors found a proof of Theorem
  \ref{thm:ultraspherical} close to the old proofs. During the preparation of the paper we found the
  references \cite{T:N}, \cite{V:L}.} Tur{\'a}n \cite{Tu}
proved that $\Delta_n^{(0)}(x)>0$ for $-1<x<1$.
The proof in \cite{T:N} of Theorem
\ref{thm:ultraspherical}  is based on a formula relating the Tur{\'a}n
determinant
$$
\Delta_{n,\lambda}(x)=F_{n,\lambda}^2(x)-F_{n-1,\lambda}(x)F_{n+1,\lambda}(x)
$$
of the normalized ultraspherical polynomials
$F_{n,\lambda}(x)=P_n^{(\lambda)}(x)/P_n^{(\lambda)}(1)$
and the expression 
$$
D_{n,\lambda}(x)=[\frac{d}{dx}P_{n}^{(\lambda)}(x)]^2-\frac{d}{dx}P_{n-1}^{(\lambda)}(x)\frac{d}{dx}P_{n+1}^{(\lambda)}(x),
$$
namely (see \cite[(5.9)]{T:N})
\begin{equation}\label{eq:TN}
\frac{\Delta_{n,\lambda}(x)}{1-x^2}=\frac{D_{n,\lambda}(x)}{n(n+2\lambda)[P_{n}^{(\lambda)}(1)]^2}.
\end{equation}
See also \cite{Da}.
Using the well-known formula for differentiation of ultraspherical polynomials
$$
\frac{d}{dx}P_n^{(\lambda)}(x)=2\lambda P_{n-1}^{(\lambda+1)}(x),
$$
we see that 
\begin{equation}
D_{n,\lambda}(x)=(2\lambda)^2\left(
[P_{n-1}^{(\lambda+1)}(x)]^2-P_{n-2}^{(\lambda+1)}(x)P_{n}^{(\lambda+1)}(x)\right).
\end{equation}
Except for the factor $(2\lambda)^2$ this is the Tur{\'a}n determinant
of order $n-1$ for the ultraspherical polynomials  corresponding to
the parameter $\lambda +1$. 

We see that this result is generalized in Theorem \ref{thm:thm1}.

 Since the proof of the monotonicity in Theorem \ref{thm:ultraspherical} depends on the fact that the
ultraspherical polynomials satisfy a differential equation, there
is little hope of extending the result to classes of orthogonal polynomials
which do not satisfy a differential equation. We have instead
attempted to find  bounds for normalized Tur{\'a}n determinants
without using monotonicity in the variable $x$.

This has also lead us to consider the following lower boundedness
condition for general orthonormal polynomials $(P_n)$:
\begin{equation}\label{eq:lowerbound}
\mbox{(LB)}\qquad \inf\{P_{n-1}^2(x)+P_n^2(x)\mid x\in\mathbb R,
n\in\mathbb N\}>0.
\end{equation}
If the condition (LB) holds then necessarily $\sum_{n=0}^\infty
P_n^2(x)=\infty$ for all $x\in\mathbb R$. Therefore the orthogonality
measure $\mu$ is uniquely determined and has no mass points.

In Proposition \ref{thm:estim} we prove  that (LB) holds for
symmetric orthonormal polynomials if the recurrence coefficients
are increasing and bounded. It turns out that for the orthonormal symmetric
Jacobi polynomials the condition (LB) holds if and only if $\alpha\ge 1/2$. 

The theory is applied to continuous $q$-ultraspherical polynomials in
Section 4. 

Concerning the general theory of orthogonal polynomials we refer the
reader to \cite{Sz},\cite{Ne},\cite{I}.

\section{Main results}

 \begin{prop}\label{thm:pro1} In addition to \eqref{eq:rec} assume
   that $\alpha_n\neq \gamma_n$ for $n=1,2\ldots.$ Then
for   $n \ge 2$ there holds $$\Delta_n =
\frac{(\gamma_n-\alpha_n)\alpha_{n-1}}{(\gamma_{n-1}-\alpha_{n-1})\gamma_n} \Delta_{n-1}\\+
 \frac{\alpha_n-\alpha_{n-1}}{(\gamma_{n-1}-\alpha_{n-1})\gamma_n}(p_{n-1}^2+p_n^2-2xp_{n-1}p_n).
 $$
 \end{prop}
\begin{proof}
By the recurrence relation we can remove either $p_{n+1}$ or $p_{n-1}$  from the formula
defining $\Delta_n.$ In this way
we obtain two equalities
\begin{eqnarray*}
\gamma_n\Delta_n&=&\alpha_{n}p_{n-1}^2+\gamma_np_n^2-xp_{n-1}p_n,\\
\alpha_n\Delta_n&=&\alpha_{n}p_n^2+\gamma_np_{n+1}^2-xp_np_{n+1}.
\end{eqnarray*}
We replace $n$ by $n-1$ in the second equality and multiply both 
 sides by $\gamma_{n }-\alpha_{n }.$
Next we subtract the resulting equality from the first one multiplied by $\gamma_{n-1 }-\alpha_{n-1 }.$
In this way we obtain after obvious simplifications
\begin{multline*}
(\gamma_{n-1 }-\alpha_{n-1 })\gamma_n\Delta_n-(\gamma_{n}-\alpha_{n})\alpha_{n-1}\Delta_{n-1}\\
=(\alpha_n\gamma_{n-1}-\alpha_{n-1}\gamma_n)(p_{n-1}^2+p_n^2)
-(\gamma_{n-1}-\gamma_n-\alpha_{n-1}+\alpha_n)\,xp_{n-1}p_n.
\end{multline*}
Taking into account that $\alpha_k+\gamma_k=1$ for $k\ge 1$ gives
$$
(\gamma_{n-1 }-\alpha_{n-1
})\gamma_n\Delta_n-(\gamma_{n}-\alpha_{n})\alpha_{n-1}\Delta_{n-1} =(\alpha_n
-\alpha_{n-1} )(p_{n-1}^2+p_n^2 -2xp_{n-1}p_n).
$$
\end{proof}

Proposition \ref{thm:pro1} implies

\begin{cor}{\bf\cite[Thm. 1]{Szw}} \label{thm:cor1} In addition to
  \eqref{eq:rec} assume that one of the following conditions holds.
\begin{enumerate}
\item[(i)] $(\alpha_n)$ is  increasing and $\alpha_n\le\gamma_n,
 \; n\ge 1$.
\item[(ii)] $(\alpha_n)$ is
 decreasing and $\alpha_n\ge \gamma_n,\;n\ge 1.$ Furthermore, assume
 that $\gamma_0=1$ or $\gamma_0\le\gamma_1/(1-\gamma_1)$.
 \end{enumerate}  Then $$\Delta_n(x)>0$$ for $-1<x<1.$
\end{cor}

\begin{proof} Assume first the additional condition
  $\alpha_n\neq\gamma_n$ for all $n\ge 0$. Since
  $p_{n-1}^2+p_n^2-2xp_{n-1}p_n> 0$ for $-1<x<1$, it suffices  in view of
Proposition \ref{thm:pro1} to show that $\Delta_1>0.$ We have
$$
\gamma_1\Delta_1(x)=\alpha_1p_0^2+\gamma_1p_1^2-xp_0p_1=
\frac{\alpha_1\gamma_0^2+(\gamma_1-\gamma_0)x^2}{\gamma_0^2},
$$
hence $\Delta_1>0$ if $\gamma_1\ge\gamma_0$.
If  $\gamma_1< \gamma_0,$ we have
$$\gamma_1\Delta_1(x)>
\gamma_1\Delta_1(1)=\frac{\alpha_1(1-\gamma_0)(\gamma_1/\alpha_1-\gamma_0)}{\gamma_0^2}.
$$
This is clearly non-negative in case (i) because $\gamma_1/\a_1\ge
1$, but also in case (ii) because of the assumptions on $\gamma_0$.

Assume next in case (i) that there is an $n$ such that
$\alpha_n=\gamma_n$ and let $n_0\ge 1$ be the smallest $n$ with this
property. Denoting $\alpha=\lim\alpha_n$, then clearly
$\alpha_n\le\alpha\le 1-\alpha\le\gamma_n$ for all $n$ and hence
$\alpha_n=\gamma_n=1/2$ for $n\ge n_0$. Therefore
$$
\Delta_n=p_{n-1}^2+p_n^2-2xp_{n-1}p_n>0
$$
for $n\ge n_0,\,-1<x<1.$ The formula of Proposition \ref{thm:pro1} can
be applied for $2\le n<n_0$ and the proof of the first case carries
over. Equality in case (ii) is treated similarly.
\end{proof}

From now on we will assume that additionally $\gamma_0=1.$ In this way the
polynomials $p_n$ are normalized at $x=1$ so that $p_n(1)=1.$ 
It follows by induction that $p_n$ has all its zeros in $]-1,1[$,
hence that the support of the orthogonality measure $\mu$ for $(p_n)$
is contained in $[-1,1]$.
Since
$p_n(-x)=(-1)^np_n(x)$ we conclude that $p_n(-1)=(-1)^n.$ Therefore for any $n\ge 0$ the
polynomial $p_{n+2}-p_n$ is divisible by $x^2-1.$ Defining
\begin{equation}\label{eq:q}
 q_n(x)=\frac{p_{n+2}(x)-p_n(x)}{x^2-1},\quad n\ge 0,
\end{equation}
$q_n(x)$ is a polynomial of degree $n.$ Moreover, an easy calculation shows that the
polynomials $q_n$ are orthogonal with respect to the probability measure
$d\nu(x)=\frac{1}{\gamma_1}(1-x^2)d\mu(x).$  By the recurrence relation (\ref{eq:rec}) and by
$\gamma_0=1$ we obtain that the polynomials $q_n$ satisfy
\begin{equation}\label{eq:qrec}
xq_n(x)=\gamma_{n+2}q_{n+1}(x)+\alpha_nq_{n-1}(x),\ n\ge 0,\, q_0=1/\gamma_1.
\end{equation}

 The following theorem contains a fundamental formula relating the
 Tur{\'a}n determinants of the polynomials
$p_n$ and $q_n.$
\begin{thm}\label{thm:thm1} For $n\ge 1$ we have
\begin{equation}\label{eq:fund}
\frac{\Delta_n(x)}{
1-x^2}=\alpha_n\gamma_nq_{n-1}^2(x)-\alpha_{n-1}\gamma_{n+1}q_{n-2}(x)q_n(x).
\end{equation}
\end{thm}
\begin{proof}
 By  \eqref{eq:rec} we get
 \begin{eqnarray*}
&& p_{k+1}-xp_k=\alpha_k(p_{k+1}-p_{k-1})=\alpha_k(x^2-1)q_{k-1},\\
&&xp_{k}-p_{k-1}=\gamma_{k}(p_{k+1}-p_{k-1})=\gamma_{k}(x^2-1)q_{k-1}.
 \end{eqnarray*}
 Therefore
 \begin{multline*}
(x^2-1)^2[\alpha_n\gamma_nq_{n-1}^2(x)-\alpha_{n-1}\gamma_{n+1}q_{n-2}(x)q_n(x)]\\
=(p_{n+1}-xp_n)(xp_{n}-p_{n-1})-(p_n-xp_{n-1})(xp_{n+1}-p_n)\\
=(1-x^2)(p_n^2-p_{n-1}p_{n+1}).
 \end{multline*}
\end{proof}

\begin{rem}\label{thm:TuranTuran} {\rm If we define
    $\tilde{q}_0=\gamma_1q_0=1$ and 
$$
\tilde{q}_n=\frac{\gamma_1\cdots\gamma_{n+1}}{\a_1\cdots\a_n}q_n,\;n\ge
1
$$
we have
\begin{equation}\label{eq:TuranTuran}
\frac{\Delta_n(x)}{1-x^2}=\frac{\gamma_n}{\a_n}\left(\frac{\a_1\cdots\a_n}{\gamma_1\cdots\gamma_n}\right)^2
\left[\tilde{q}_{n-1}^2(x)-\tilde{q}_{n-2}(x)\tilde{q}_n(x)\right],
\end{equation}
showing that the normalized Tur{\'a}n determinant \eqref{eq:fund} is proportional to
a Tur{\'a}n determinant of order $n-1$ of the renormalized polynomials
$(\tilde{q}_n)$. They satisfy the recursion equation ($\tilde{q}_{-1}:=0$)
$$
x\tilde{q}_n=\a_{n+1}\tilde{q}_{n+1}+\gamma_{n+1}\tilde{q}_{n-1},\;n\ge 0.
$$ 
}
\end{rem}

\begin{thm}\label{thm:thm2}
Assume that $(p_n)$ satisfies \eqref{eq:rec} with $\gamma_0=1.$ Let $(\alpha_n)$ be increasing,
$\alpha_n\le 1/2$ and
\begin{equation}\label{eq:qconcave} \alpha_n-\alpha_{n-1}\ge \frac{\alpha_n}{
1-\alpha_n}(\alpha_{n+1}-\alpha_n),\ n\ge 1.
\end{equation} (e.g. (\ref{eq:qconcave}) is satisfied if $\alpha_n$ is
 concave). Then for
$\Delta_n(x)$ defined by \eqref{eq:Turan} we have
$$
 \frac{\Delta_n(x)}{1-x^2}\ge c\Delta_n(0),\qquad -1< x< 1,\ n\ge 1,
$$
where $c={2\alpha_1\gamma_2/ \gamma_1}.$
\end{thm}
\begin{proof} Observe that (\ref{eq:qconcave}) is equivalent to 
$(\alpha_n\gamma_{n+1})$ being increasing. Let
$$D_n(x)= \gamma_nq_{n-1}^2(x)- \gamma_{n+1}q_{n-2}(x)q_n(x).$$
Since $\alpha_n\ge \alpha_{n-1}$ Theorem \ref{thm:thm1}  implies that
\begin{equation}\label{eq:aux}
\frac{\Delta_n(x)}{1-x^2}\ge \alpha_{n-1}D_{n}(x).
 \end{equation}
 By \eqref{eq:qrec} we can remove $q_{n}$ or
$q_{n-2}$ from the expression  defining $D_n.$ In this way we obtain
\begin{eqnarray}
D_n&=&\alpha_{n-1}q_{n-2}^2+\gamma_nq_{n-1}^2-xq_{n-2}q_{n-1},\label{eq:D}\\
\frac{\alpha_{n-1}}{\gamma_{n+1}}D_n&=&\frac{\alpha_{n-1}\gamma_n}{\gamma_{n+1}}q_{n-1}^2+
\gamma_{n+1}q_n^2-xq_{n-1}q_n.\nonumber
\end{eqnarray}
Replacing $n$ by $n-1$ in the second equality and subtracting it from
the first we find
$$
D_n-\frac{\alpha_{n-2}}{\gamma_{n}}D_{n-1}=
\frac{\alpha_{n-1}\gamma_n-\alpha_{n-2}\gamma_{n-1}}{\gamma_n}q_{n-2}^2\ge
0.
$$
By iterating the  inequality between $D_n$ and $D_{n-1}$ we obtain
 \begin{equation}\label{eq:use}
D_n\ge
{\alpha_1\ldots\alpha_{n-2}\over \gamma_3\ldots\gamma_n}D_2.
\end{equation}
From  \eqref{eq:D} we get 
\begin{equation}\label{eq:use2}
D_2=\alpha_1q_0^2+\gamma_2q_1^2-xq_0q_1
={\alpha_1\over \gamma_1^2}+{\gamma_2x^2\over
\gamma_1^2\gamma_2^2}-{x^2\over\gamma_1^2\gamma_2} ={\alpha_1\over
\gamma_1^2},
\end{equation} 
so \eqref{eq:aux} implies
$$ {\Delta_n(x)\over
1-x^2}\ge {\alpha_1\ldots\alpha_{n}\over \gamma_1\ldots\gamma_n}{\alpha_1\gamma_2\over
\alpha_n\gamma_1}\ge {\alpha_1\ldots\alpha_{n}\over \gamma_1\ldots\gamma_n}{2\alpha_1\gamma_2\over
 \gamma_1}.
$$
 The conclusion follows from the next lemma.\end{proof}
\begin{lem}\label{thm:deltabounds}
Under the assumptions of Theorem \ref{thm:thm2} there holds
$$
\Delta_n(0)\le{\alpha_1\ldots\alpha_{n }\over \gamma_1\ldots\gamma_n}\le {\gamma_1\over\alpha_1}
\Delta_n(0)  ,\quad n\ge 1.
$$
\end{lem}
\begin{proof}
Denote $$h_n={\gamma_1\ldots\gamma_n\over \alpha_1\ldots\alpha_{n }}.$$  By (\ref{eq:rec})
we have
$$p_{2n}(0)=(-1)^n{\alpha_1\alpha_3\ldots\alpha_{2n-1}\over
\gamma_1\gamma_3\ldots\gamma_{2n-1}}.$$ Hence
$$
\Delta_{2n}(0)h_{2n}=p_{2n}^2(0)h_{2n}=  \prod_{k=1}^n{\alpha_{2k-1}\over\alpha_{2k}} \prod_{k=1}^{n } {\gamma_{2k
}\over\gamma_{2k-1}}\le 1.
$$
 On the other hand
\begin{multline*} \Delta_{2n+1}(0)h_{2n+1}= -p_{2n}(0)p_{2n+2}(0)h_{2n+1}\\=
\prod_{k=1}^n{\alpha_{2k-1}\over \alpha_{2k}}\prod_{k=1}^n{\gamma_{2k}\over
\gamma_{2k-1}} =\Delta_{2n}(0)h_{2n}\le 1.
\end{multline*}
Moreover
\begin{multline*}
\Delta_{2n}(0)h_{2n}=\prod_{k=1}^n{\alpha_{2k-1}\over\alpha_{2k}} \prod_{k=1}^{n }
{\gamma_{2k }\over\gamma_{2k-1}} \ge \prod_{k=2}^{2n}{\alpha_{k-1}\over\alpha_{k}}
\prod_{k=2}^{2n } {\gamma_{k
}\over\gamma_{k-1}}={\alpha_1\gamma_{2n}\over\gamma_1\alpha_{2n}}\ge {\alpha_1\over
\gamma_1}.
\end{multline*}
\end{proof}

Theorem \ref{thm:thm2} has the following counterpart and the proof
is very similar: 

\begin{thm}\label{thm:thm2a}
Assume that $(p_n)$ satisfies \eqref{eq:rec} with $\gamma_0=1.$ Let
$\alpha_n,n\ge 1$ be decreasing,
$\alpha_n\ge {1\over 2}$ and
\begin{equation}\label{eq:qconcave2}
 \alpha_n-\alpha_{n-1}\le {\alpha_n\over
1-\alpha_n}(\alpha_{n+1}-\alpha_n),\ n\ge 2.
\end{equation}
Then for
$\Delta_n(x)$ defined by \eqref{eq:Turan} we have
$$
 {\Delta_n(x)\over 1-x^2}\le C\Delta_n(0),\qquad -1< x< 1,\ n\ge 1,
$$
where $C=2\gamma_2.$ (Note that \eqref{eq:qconcave2} implies convexity
of $\a_n,n\ge 1$.)
\end{thm}

\begin{rem}\label{thm:2}{\rm

Note that the normalized symmetric Jacobi polynomials
$p_n(x)=R_n^{(\a,\a)}(x)$ given by \eqref{eq:Jacobinorm} satisfy \eqref{eq:rec} with
\begin{equation}\label{eq:recsymJacobi}
\gamma_n=\frac{n+2\a+1}{2n+2\a+1},\quad \a_n=\frac{n}{2n+2\a+1}.
\end{equation}
(In the case of $\a=-1/2$, i.e. Chebyshev polynomials of the first
kind, these formulas shall be interpreted as $\gamma_0=1,\a_0=0$.)

 For $\alpha\ge -1/2$ we have $(\a_n)$ is increasing and concave and
$c=1$. 

For $-1<\a\le -1/2$ the sequence $(\a_n)$ is decreasing,
\eqref{eq:qconcave2} holds and $C=1$.

 The statement about the
constants $c,C$ follows from
Corollary \ref{thm:bounds}.  However we cannot expect $c=1$ in
general, because it is easy to construct an example, where the
normalized Tur{\'a}n determinant \eqref{eq:fund} is not monotone for $0<x<1$.

Consider the sequence
$(\a_n)=(0,1/2-3\varepsilon,1/2-2\varepsilon,1/2-\varepsilon,1/2,1/2,\ldots)$,
which is increasing and concave for $0<\varepsilon<1/8$.
In this case the Tur{\'a}n determinant
$\tilde{q}_2^2-\tilde{q}_1\tilde{q}_3$ is proportional to
$f(x)=x^4+A(\varepsilon)x^2 +B(\varepsilon)$, where
$$
A(\varepsilon)=4\varepsilon^2+3\varepsilon-1/2,\; B(\varepsilon)=(1/2-3\varepsilon)^2(1/2-\varepsilon)(1/2+2\varepsilon)^2/\varepsilon.
$$
Clearly, $f$ is not monotone for $0<x<1$ when $\varepsilon$ is small.

}
\end{rem}

\begin{cor}\label{thm:density} Under the assumptions of Theorem
  \ref{thm:thm2} and the additional hypothesis $\lim\,\a_n=1/2$, the
  orthogonality measure $\mu$ is
absolutely continuous on $(-1,1)$ with a strictly positive and continuous density $g(x)=d\mu(x)/dx$ satisfying
$$
g(x)\le {C\over \sqrt{1-x^2}}.
$$
\end{cor}

\begin{proof} The corresponding orthonormal polynomials $(P_n)$
  satisfy
\begin{equation}\label{eq:ONP}
xP_n=\lambda_nP_{n+1}+\lambda_{n-1}P_{n-1},
\end{equation}
where $\lambda_n=\sqrt{\alpha_{n+1}\gamma_n}$. We also have
$P_n=\delta_np_n$, where
$$
\delta_n=\sqrt{\frac{\gamma_0\cdots
    \gamma_{n-1}}{\a_1\cdots\a_n}},\;n\ge 1,\;\;\delta_0=1,
$$
 and $\lim\lambda_n=1/2$.
Since
$$
\lambda_{n+1}-\lambda_n=
\frac{\a_{n+2}(\gamma_{n+1}-\gamma_n)+\gamma_n(\a_{n+2}-\a_{n+1})}
{\sqrt{\a_{n+2}\gamma_{n+1}}+\sqrt{\a_{n+1}\gamma_{n}}},
$$
the monotonicity of $(\a_n),(\gamma_n)$ implies that
 $$
\sum_{n=1}^\infty|\lambda_{n+1}-\lambda_n|<\infty.
$$
By the theorem in  \cite{M:N} we conclude that the orthogonality
measure $\mu$ has a positive continuous density $g(x)$ for
$-1<x<1$. Furthermore, it is known from
this theorem that
$$
\lim_{n\to\infty} [P_n^2(x)-P_{n-1}(x)P_{n+1}(x)]=
\frac{2\sqrt{1-x^2}}{\pi g(x)},
$$
unifomly on compact subsets of $]-1,1[$.
For another proof of this result see \cite[p. 201]{F:L:S}, where it is
also proved that $(P_n(x))$ is uniformly bounded on compact subsets of
$]-1,1[$ for $n\to\infty$.
We have
$$
\Delta_n(x)=\frac{1}{\delta_n^2}\left(P_n^2(x)-k_nP_{n-1}(x)P_{n+1}(x)\right),
$$
where
$$
k_n=\frac{\delta_n^2}
{\delta_{n-1}\delta_{n+1}}=\sqrt{\frac{\a_{n+1}\gamma_{n-1}}{\a_n\gamma_n}},
$$
and it follows that $\lim k_n=1$. Using that
$$
\frac{\Delta_n(x)}{\Delta_n(0)}=\frac{P_n^2(x)-k_nP_{n-1}(x)P_{n+1}(x)}{P_n^2(0)-k_nP_{n-1}(0)P_{n+1}(0)},
$$
we get the result.
\end{proof}

In analogy with the proof of Corollary \ref{thm:density} we get

\begin{cor}\label{thm:densityrev} Under the assumptions of Theorem
  \ref{thm:thm2a} and the additional hypothesis $\lim\,\a_n=1/2$, the
  orthogonality measure $\mu$ is
absolutely continuous on $(-1,1)$ with a strictly positive and continuous density $g(x)=d\mu(x)/dx$ satisfying
$$
g(x)\ge {C\over \sqrt{1-x^2}}.
$$
\end{cor}

We now return to the polynomials $(q_n)$ defined in \eqref{eq:q} and
prove that they have a non-negative Tur{\'a}n determinant after
normalization to being 1 at 1.
The polynomials $q_n$ are orthogonal with respect to a measure
supported by $[-1,1].$ Therefore $q_n(1)>0.$

\begin{prop}\label{thm:Turanq} Under the assumptions of Theorem
  \ref{thm:thm2} we have for $n\ge 1$
$${q_n^2(x)\over q_n^2(1)}-{q_{n-1}(x)\over q_{n-1}(1)}{q_{n+1}(x)\over q_{n+1}(1)}\ge
0.$$
\end{prop}
\begin{proof}
   Indeed, let $Q_n(x)=\displaystyle{q_n(x)\over q_n(1)}.$ Then
$$xQ_n=c_nQ_{n+1}+(1-c_n)Q_{n-1},$$ where
$$c_n=\gamma_{n+2}{q_{n+1}(1)\over q_n(1)}.$$
We will  show that $c_n$ is decreasing and $c_n\ge 1/2.$ Then the conclusion follows
from Corollary \ref{thm:cor1}. But $c_{n-1}\ge c_n$ is equivalent to
$$
D_{n+1}(1)=\gamma_{n+1}q_n^2(1)- \gamma_{n+2} q_{n-1}(1)q_{n+1}(1)\ge
0,
$$
which follows from \eqref{eq:use} and \eqref{eq:use2}. We will show that $c_n\ge {1/2}$ by
induction. We have
$$
c_0=\gamma_2{q_1(1)\over q_0(1)}=1.
$$
Assume $c_{n-1}\ge 1/2.$ 
By \eqref{eq:qconcave} the sequence $(\alpha_n\gamma_{n+1})$ is increasing. Putting
$\alpha=\lim \alpha_n$ we then get 
$$
\alpha_n\gamma_{n+1}\le \alpha(1-\alpha)\le \frac14.
$$
Using this and \eqref{eq:qrec} leads to
$$
1= c_n+{\alpha_n\gamma_{n+1}\over c_{n-1}}\le c_n+{1\over 4c_{n-1}}\le
c_n+{1\over 2},
$$
hence $c_n\ge 1/2.$
\end{proof}

\section{Lower bound estimates}
It turns out that Tur\'an determinants can be used to obtain lower bound estimates for
orthonormal polynomials. Recall that if the polynomials $p_n$ satisfy the recurrence
relation \eqref{eq:rec}, then their orthonormal version $(P_n)$ satisfy
$$
xP_n=\lambda_nP_{n+1}+\lambda_{n-1}P_{n-1},
$$
where $\lambda_n=\sqrt{\alpha_{n+1}\gamma_n}.$

\begin{prop}\label{thm:estim}
Assume that the  polynomials $(P_n(x))$ satisfy
\begin{equation}\label{eq:recorth}
xP_n=\lambda_nP_{n+1}+\lambda_{n-1}P_{n-1},\ n\ge 0,
\end{equation}
with $P_{-1}=\lambda_{-1}=0,$ $\lambda_n>0,\ n\ge 0,$ and $P_0=1.$ If the sequence
$(\lambda_n)$ is increasing and $\lim\lambda_n=L<\infty,$ then the
(LB) condition \eqref{eq:lowerbound} holds, viz.
$$P_n^2(x)+P_{n-1}^2(x)\ge \frac{\lambda_0^2}{2L^2}.
 $$
\end{prop}
\begin{proof}
This proof is inspired by \cite[Thm. 3]{askey}. By replacing the
polynomials $P_n(x)$ by $P_n(2Lx)$ we can assume that $\lim
\lambda_n=1/2$. This assumption
implies that the corresponding Jacobi
matrix is a contraction, because it can be majorized by the Jacobi matrix with entries
$\lambda_n={1\over 2}.$ Therefore the orthogonality measure is supported by the interval
$[-1,1].$ In this way it suffices to consider $x$ from $[-1,1]$ because the functions
$P_n^2(x)$ are increasing on $[1,+\infty[$ and $P_n^2(-x)=P_n^2(x).$ Let
$$
\mathcal D_n(x)=\lambda_{n-1}P_n^2(x)-\lambda_nP_{n-1}(x)P_{n+1}(x),\quad n\ge 1.$$
By  \eqref{eq:recorth} we can remove $P_{n+1}$ to get
\begin{equation}\label{eq:one}
\mathcal D_n=\lambda_{n-1}P_{n-1}^2+\lambda_{n-1}P_n^2-xP_{n-1}P_{n}.\end{equation}
 Alternatively we can remove
$P_{n-1}$  and obtain
\begin{equation}\label{eq:two}
 \frac{\lambda_{n-1}}{\lambda_n}\mathcal D_n=\lambda_nP_{n+1}^2+
{\lambda_{n-1}^2\over
\lambda_n}P_{n}^2-xP_nP_{n+1}.
\end{equation}
Replacing $n$ by $n-1$ in \eqref{eq:two} and subtracting it from (\ref{eq:one}) gives
\begin{equation}\label{eq:three}
\mathcal  D_n-  {\lambda_{n-2}\over \lambda_{n-1}}\mathcal D_{n-1}={\lambda_{n-1}^2-\lambda_{n-2}^2\over
\lambda_{n-1}}P_{n-1}^2\ge 0.
\end{equation}
By iterating the inequality
$\mathcal D_n\ge(\lambda_{n-2}/\lambda_{n-1})\mathcal D_{n-1}$,
 we obtain
$$
\mathcal D_n\ge {\lambda_0\over \lambda_{n-1}}\mathcal
D_1={\lambda_0^2\over \lambda_{n-1}}\ge 2\lambda_0^2,
$$
because by (\ref{eq:one}) we have $\mathcal D_1=\lambda_0.$ Now  (\ref{eq:one}) implies for $|x|\le 1$
\begin{equation}\label{eq:help}
\mathcal D_n\le \lambda_{n-1}P_{n-1}^2+\lambda_{n-1}P_n^2+ {1\over 2}|x|(P_{n-1}^2+P_n^2)\le
P_{n-1}^2+P_n^2.
\end{equation}
In the general case the lower bound is $2(\lambda_0/(2L))^2$.
\end{proof}

\begin{cor}\label{thm:density2} Under the assumptions of Proposition
  \ref{thm:estim} with $L=1/2$ the orthogonality measure $\mu$ is
absolutely continuous  with a continuous density $g=d\mu(x)/dx$ on
$[-1,1]$ satisfying
$$
g(x)\le {{1\over 2\pi\lambda_0^2}  \sqrt{1-x^2}}.
$$
Furthermore, $g(x)>0$ for $-1<x<1$.
\end{cor}

\begin{proof}
By assumptions the orthogonality measure is supported by $[-1,1]$.
 By the proof of Proposition \ref{thm:estim} we have
$$
\mathcal D_n(x)\ge 2\lambda_0^2.
$$
On the other hand by \cite{M:N} and \cite[p. 201]{F:L:S} the orthogonality measure is absolutely continuous in
the interval $]-1,1[$ with a strictly positive and continuous density
$g$ such that
$$
 \lim_{n\to\infty}{1\over \lambda_{n-1}} \mathcal D_n(x) =
 {2\sqrt{1-x^2}\over \pi g(x)},
$$
uniformly on compact subsets of $]-1,1[$, cf. the proof of Corollary \ref{thm:density}.
By Property (LB) there are no masses at $\pm 1$.
\end{proof}

The Jacobi polynomials $P_n^{(\a,\a)}(x)$ in the standard notation
 of Szeg{\H o}, cf. \cite{Sz}, are discussed in the Introduction.
 The corresponding orthonormal polynomials are
 denoted $P_n(\a;x)$. We recall that
\begin{equation}\label{eq:JacobiL2}
c_\a\int_{-1}^1
[P_n^{(\a,\a)}(x)]^2(1-x^2)^{\a}\,dx=\frac{2^{2\a+1}\Gamma(n+\a+1)^2}
{(2n+2\a+1)n!\Gamma(n+2\a+1)B(\a+1,1/2)}.
\end{equation}

\begin{rem}\label{thm:DomNevai} {\rm Corollary \ref{thm:density2} is
    also obtained in \cite[p.758]{D:N}. 
}
\end{rem}

\begin{cor}\label{thm:ultrabd} Let $(P_n(\a;x))$ denote the orthonormal
  symmetric Jacobi polynomials.
\begin{enumerate}
\item[(i)] For $\alpha\ge 1/2$ we have
$$\inf\{P_n^2(\a;x)+P_{n-1}^2(\a;x)\ |\ x\in \mathbb{R},n\in\mathbb{N}\}\ge {2\over 2\alpha+3}  .$$
\item[(ii)] For $-1<\alpha<1/2$ we have
$$\inf\{P_n^2(\a;x)+P_{n-1}^2(\a;x)\ |\ x\in \mathbb{R},n\in\mathbb{N}\}=0.$$
\end{enumerate}
\end{cor}
\begin{proof}
Assume $\alpha\ge 1/2.$ In this case we get from \eqref{eq:recsymJacobi}
$$
\lambda_n^2={1\over 4} \left [1-{4\alpha^2-1\over
    4(n+\alpha+1)^2-1}\right ],
$$
so $(\lambda_n)$ is increasing with $\lim\,\lambda_n=1/2$. By Proposition \ref{thm:estim} we thus have
$$P_n^2+P_{n-1}^2\ge 2\lambda_0^2={2\over 2\alpha+3},$$
which shows (i).

In order to show (ii) we will make use of Hilb's asymptotic formula
\cite[Thm 8.21.12]{Sz}:

\begin{multline}\label{Hilb}\theta^{-1/2}\left (\sin{\theta\over 2}\right )^{\alpha+1/2} \left (\cos{\theta\over
2}\right )^{\alpha+1/2}P_n^{(\alpha,\alpha)}(\cos\theta)\\= {\Gamma(\alpha+n+1)\over n!
\sqrt{2}N^\alpha} J_\alpha(N\theta) +O(n^{-3/2}),
\end{multline}
 where $\theta\in [c/n,\pi/2]$,
$N=n+\alpha+{1\over 2}$ and $c>0$ is fixed.  Let $j_\a$ denote the smallest positive zero of the
Bessel function $J_\alpha$.

 Defining $\theta_n=j_\a/N$ we get
\begin{eqnarray*}
n^{-\a}P_n^{(\alpha,\alpha)}(\theta_n)&=& O(n^{-3/2}),\\
n^{-\a}P_{n-1}^{(\alpha,\alpha)}(\theta_n)&=&(1/\sqrt{2}+o(1))J_\a(j_\a\frac{n+\a-1/2}{n+\a+1/2})+ O(n^{-3/2})=O(n^{-1}).
\end{eqnarray*}
By \eqref{eq:JacobiL2} and Stirling's formula 
$$
c_\a\int_{-1}^1[P_n^{(\alpha,\alpha)}(x)]^2(1-x^2)^\a\,dx\sim \frac{2^{2\a+1}}{B(\a+1,1/2)}n^{-1},
$$
and hence
$$
P_n^2(\a;\cos\theta_n)=O(n^{2\a-2}),\quad
P_{n-1}^2(\a;\cos\theta_n)=O(n^{2\a-1}).
$$
This shows that
$$
P_n^2(\a;\cos\theta_n)+P_{n-1}^2(\a;\cos\theta_n)\to 0 \mbox{ when }
\a<1/2.
$$
\end{proof}
\begin{rem}\label{thm:Rem3} {\rm The example of symmetric Jacobi polynomials suggests that if $\lambda_n$ is
decreasing, then the condition of Corollary \ref{thm:ultrabd} (ii) may hold. This is not true
because for $\frac12<\lambda_0 < {1\over\sqrt 2}$ and $\lambda_n={1\over
  2}$ for $n\ge 1$ we have a decreasing sequence and the
corresponding Jacobi matrix has norm 1 because this is so for the cases
$\lambda_0=\frac12$ and $\lambda_0=1/\sqrt{2}$, which correspond to the
Chebyshev polynomials of the second and first kind
respectively. Furthermore, for $n\ge 2$ we have by \eqref{eq:one}  and \eqref{eq:three}
$$
\mathcal D_n=
  \lambda_{n-1}P_n^2-{\lambda_n}P_{n-1}P_{n+1}=\mathcal D_2={2\over \lambda_0^2}
[\lambda_0^4-(\lambda_0^2-\frac14)x^2]
$$
and for $-1<x<1$
$$
\mathcal D_2(x)>\mathcal D_2(1)={2\over
  \lambda_0^2}(\lambda_0^2-\frac12)^2>0.
$$
On the other hand \eqref{eq:help} applies for $n\ge 2$ and we see that
the orthonormal polynomials satisfy
$$
\inf\{P_n^2(x)+P_{n-1}^2(x)\mid x\in\mathbb R,n\in\mathbb N\}\ge {2\over
  \lambda_0^2}(\lambda_0^2-\frac12)^2.
$$
}
\end{rem}

\section{Continuous $q$-ultraspherical polynomials}

The continuous $q$-ultraspherical polynomials $C_n(x;\b|q)$ depend on
two real parameters $q,\b$, and for $|q|,|\b|<1$ they are orthogonal
with respect to a continuous weight function on $]-1,1[$,
cf. \cite{I},\cite{K:S}. The 3-term recurrence relation is
\begin{equation}\label{eq:qultra}
xC_n(x;\b|q)=\frac{1-q^{n+1}}{2(1-\b q^n)}C_{n+1}(x;\b|q)+\frac{1-\b^2
  q^{n-1}}{2(1-\b q^n)}C_{n-1}(x;\b|q),\quad n\ge 0
\end{equation}
with $C_{-1}=0,C_0=1$. The orthonormal version $\mathcal C_n(x;\b|q)$
satisfy equation \eqref{eq:recorth} with
\begin{equation}\label{eq:qultraon}
\lambda_n=\frac{1}{2}\sqrt{\frac{(1-q^{n+1})(1-\b^2 q^n)}{(1-\b
  q^n)(1-\b q^{n+1})}}.
\end{equation}
The value $C_n(1;\b|q)$ is not explicitly known, and therefore we can
only obtain the recurrence coefficients $\a_n,\gamma_n$ from
\eqref{eq:rec} for
$p_n(x)=C_n(x;\b|q)/C_n(1;\b|q)$ as given by the recursive equations
\begin{equation}\label{eq:qultra1}
\a_{n+1}=\frac{\lambda_n^2}{1-\a_n},\quad \a_0=0,\; \gamma_n=1-\a_n,
\end{equation}
which we get from the relation $\lambda_n=\sqrt{\a_{n+1}\gamma_n}$.

\begin{thm}\label{thm:qultra} (i) Assume $0\le\b\le q<1$. Then the recurrence coefficients
  $(\lambda_n)$ form an increasing sequence with limit $1/2$, and
  therefore $(\mathcal C_n(x;\b|q))$ satisfies (LB).

(ii) Assume $0\le q\le\b<1$. Then the recurrence coefficients
  $(\lambda_n)$ form a  decreasing sequence with limit $1/2$, and the sequence
  $(\a_n)$ is increasing and concave with limit $1/2$. In particular,
we have
$$
\frac{\Delta_n(x)}{1-x^2}\ge c\Delta_n(0),\quad -1<x<1, n\ge 1,
$$
with $c=2\a_1(1-\a_2)/(1-\a_1)$.
\end{thm}

\begin{proof} The function
$$
\psi(x)=\frac{(1-qx)(1-\b^2x)}{(1-\b x)(1-\b
  qx)}=1+(1-\b)(\b-q)\frac{x}{(1-\b x)(1-\b qx)}
$$
is decreasing for $0\le\b\le q<1$ and increasing for $0\le q\le
\b<1$. This shows that $\lambda_n=(1/2)\sqrt{\psi(q^n)}$ is increasing
in case (i) and decreasing in case (ii). In both cases the limit is
$1/2$.

In case (ii) we therefore have $\lambda_n^2\ge 1/4$ and hence 
$$
\a_{n+1}\ge\frac{1}{4(1-\a_n)}\ge \a_n,
$$
because $4x(1-x)\le 1$ for $0\le x\le 1$. This shows that $(\a_n)$ is
increasing and hence with limit $1/2$. We further have
$$
\a_{n+1}-\a_n=2(\lambda_n^2-\frac{1}{4})+2(\frac12-\a_n)(\frac12-\a_{n+1}),
$$
which shows that $\a_{n+1}-\a_n$ is decreasing, i.e. $(\a_n)$ is concave.
We can now apply Theorem \ref{thm:thm2}.
\end{proof}

\noindent
Christian Berg\\
Department of Mathematics, University of Copenhagen,
Universitetsparken 5, DK-2100, Denmark\\
e-mail: {\tt{berg@math.ku.dk}}

\vspace{0.4cm} \noindent
Ryszard Szwarc\\
Institute of Mathematics, University of Wroc{\l}aw, pl.\ Grunwaldzki 2/4, 50-384
Wroc{\l}aw, Poland
\newline  and
\newline \noindent Institute of Mathematics and Computer Science,
University of Opole, ul. Oleska 48,
45-052 Opole, Poland\\
e-mail: {\tt{szwarc2@gmail.com}}

 \end{document}